\documentclass[12pt,a4paper,fleqn]{article}
\usepackage{amsfonts,amsmath,latexsym,amssymb,graphicx,mathrsfs}

\usepackage{graphicx}
\usepackage{color}
\usepackage{amssymb}
\usepackage{amssymb}
\usepackage[T1]{fontenc}
\usepackage{latexsym}
\usepackage{xypic}
\usepackage{eufrak}
\usepackage{euscript}
\usepackage{amsfonts,amsmath}
\usepackage{verbatim}
\usepackage[english]{babel}
\usepackage{mathrsfs}

\newtheorem{prop}{Proposition}[section]

\newtheorem{lemma}[prop]{Lemma}
\newtheorem{rem}[prop]{Remark}

\newtheorem{thm}[prop]{Theorem}

\renewcommand{\geq}{\geqslant}
\def\leq{\leqslant}

\newcommand{\Z}{\mathbb{Z}}
\newcommand{\R}{\mathbb{R}}

\newcommand{\C}{\mathbb{C}}

\def\HH{\EuFrak H}

\def\e{\varepsilon}

\def\1{{\mathbf{1}}}

\def\1{{\mathbf{1}}}
\def\0.5{{\frac{1}{2}}}

\newcommand{\beas}{\begin{eqnarray*}}
\newcommand{\enas}{\end{eqnarray*}}

\newcommand{\NN}{\mathscr{N}}

\newcommand{\fin}
{ \vspace{-0.6cm}
\begin{flushright}
\mbox{$\Box$}
\end{flushright}
\noindent }

\newcommand{\qed}{\nopagebreak\hspace*{\fill}
{\vrule width6pt height6ptdepth0pt}\par}

\begin{document}

\begin{center}
{\Large{\bf A multiple stochastic integral criterion \\ for almost sure limit theorems}}
\\~\\ by 
Bernard Bercu\footnote{Institut de Math\'ematiques de Bordeaux, Universit\'e Bordeaux 1, 351 cours de
la lib\'eration, 33405 Talence cedex, 
France. Email: \texttt{Bernard.Bercu@math.u-bordeaux1.fr}},
Ivan Nourdin\footnote{Laboratoire de Probabilit\'es et Mod\`eles Al\'eatoires, Universit\'e
 Pierre et Marie Curie (Paris VI), Bo\^ite courrier 188, 4 place Jussieu, 75252 Paris Cedex 05, 
France. Email: \texttt{ivan.nourdin@upmc.fr}}
and Murad S. Taqqu\footnote{Boston University, Departement of Mathematics, 111 Cummington Road, Boston (MA), USA. 
Email: \texttt{murad@math.bu.edu}}\footnote{Murad S. Taqqu was partially supported by the NSF Grant DMS-0706786 at
Boston University.}\\ 
{\small {\it Universit\'e Bordeaux 1, Universit\'e Paris 6 and Boston University}}\\~\\
\end{center}

{\small \noindent {\bf Abstract:} 
In this paper, we study almost sure central limit theorems for multiple stochastic integrals and provide
a criterion based on the kernel of these multiple integrals. We apply our result to normalized partial sums
of Hermite polynomials of increments of fractional Brownian motion. We obtain almost sure central limit
theorems for these normalized sums when they converge in law to a normal distribution.\\
}

{\small \noindent {\bf Key words}: Almost sure limit theorem; multiple stochastic integrals;
fractional Brownian motion; Hermite power variation.
\\


\noindent {\it This version}: April 14, 2009}

\section{Introduction}
Let $\{X_{n}\}_{n\geq 1}$ be a sequence of real-valued independent identically distributed random variables with 
$E[X_{n}]=0$ and $E[X_{n}^{2}]=1$, and denote
\[
S_{n}= \frac1{\sqrt{n}}\sum_{k=1}^n  X_k.
\]
The celebrated almost sure central limit theorem (ASCLT) 
states that the sequence of random empirical measures, given by 
\[
\frac{1}{\log n}\sum_{k=1}^{n}\frac{1}{k}\delta_{S_k},
\]
converges almost surely to the $\NN(0,1)$ distribution as $n\to\infty$.
In other words, if $N$ is a $\NN(0,1)$ random variable, then, almost surely, for all $x\in\R$,
\[
\frac{1}{\log n}\sum_{k=1}^{n}\frac{1}{k}{\bf 1}_{\{S_k\leq x\}} \longrightarrow 
P(N\leq x),\quad\mbox{as $n\to\infty$}.
\]
or, equivalently, almost surely, for any bounded and continuous function $\varphi:\R\to\R$, 
\begin{equation}
\label{ASCLTh}
\frac{1}{\log n}\sum_{k=1}^{n}\frac{1}{k}\varphi(S_{k}) \longrightarrow 
E[\varphi(N)],\quad\mbox{as $n\to\infty$}.
\end{equation}

The ASCLT was stated first by L\'evy \cite{Levy} without proof.
It was then forgotten for half century.
It was rediscovered by Brosamler \cite{Brosamler} and Schatte \cite{Schatte}
and proven, in its present form, by Lacey and Philipp \cite{LaceyPhillip}. 
We refer the reader to Berkes and Cs\'aki \cite{BerkesCsaki} for a universal
ASCLT covering a large class of limit theorems for partial sums, extremes, empirical
distribution function and local times associated with independent random variables $\{X_n\}$,
as well as the PhD thesis by
Gonchigdanzan \cite{these}, where extensions of the ASCLT to weakly dependent
random variables are studied, for example in the context of strong mixing or $\rho$-mixing.
Ibragimov and Lifshits \cite{IL,
IbragimovLifshits2} have provided a criterion for (\ref{ASCLTh}) 
which does not require the sequence $\{X_n\}$ of random variables to be necessarily independent
nor that the sequence $\{S_n\}$ have the specific form of partial sums.
This criterion is stated in Proposition \ref{thm-IL} below.

Our goal is to investigate the ASCLT for a sequence of multiple stochastic integrals. Conditions ensuring
the convergence in law of this sequence to the standard $\NN(0,1)$ distribution are now well-known, see
Nuarlart and Peccati \cite{NP}. We will derive a criterion for this sequence of multiple
integrals to satisfy also the ASCLT.

As an application, we consider some non-linear
functions of strongly dependent Gaussian random variables. 
We will generate strong dependence by using increments of a standard fractional Brownian motion $B^H$.   
Recall that $B^H=(B^H_t)_{t\geq 0}$ is a centered Gaussian process with continuous paths such that 
\[
E[B^H_tB^H_s]=\frac{1}{2}\Bigl(t^{2H}+s^{2H} -|t-s|^{2H}\Bigr),\quad s,t\geq 0.
\]
The process $B^H$ is self-similar with stationary increments
and we refer the reader to Nualart \cite{Nbook} and Samorodnitsky and Taqqu \cite{SamorTaqqu}
for its main properties. The increments 
\[
Y_k=B^H_{k+1}-B^H_{k},\quad k\geq 0,
\]
called ``fractional Gaussian noise'',
are centered stationary Gaussian random variables with covariance
\begin{equation}\label{e:cov}
\rho(r)=E[Y_kY_{k+r}]=\frac12\big(|r+1|^{2H}+|r-1|^{2H}-2|r|^{2H}\big), \quad r\in\mathbb{Z}.
\end{equation} 
This covariance behaves asymptotically as
$$
\rho(r)\sim H(2H-1)|r|^{2H-2}\quad\mbox{as $|r|\to\infty$}.
$$
Observe that $\rho(0)=1$ and $(i)$ for $0<H<1/2$: $\rho(r)<0$ for $r\neq 0$, $\sum_{r\in\Z}|\rho(r)|<\infty$ and 
$\sum_{r\in\Z}\rho(r)=0$; $(ii)$ for $H=1/2$: $\rho(r)=0$ if $r\neq 0$; and $(iii)$
for $1/2<H<1$: $\sum_{r\in\Z}|\rho(r)|=\infty$.
The Hurst index measures the strenght of the dependence when $H\geq 1/2$: the larger $H$,
the stronger the dependence.  
%Berkes and Horv\'ath \cite{BerkesHorvath} prove an ASCLT
%for fractional Brownian motion. They show that, for any bounded continuous real function $\varphi$,
%\begin{equation}\label{berkes}
%\frac{1}{\log T}\int_{1}^{T}\frac{1}{t}\varphi\Bigl(\frac{B^H_{t}}{t^H}\Bigr) \overset{\rm a.s}{\longrightarrow} 
%E[\varphi(N)]
%\quad\mbox{as $T\to\infty$},
%\end{equation}
%with $N\sim\NN(0,1)$.
%They also obtain a sharp strong approximation result for this continuous time version of the ASCLT. 

We shall consider random variables
\begin{equation}\label{e:XY}
H_q(Y_k)=H_q(B^H_{k+1}-B^H_k),\quad k\geq 0,
\end{equation}
where $H_q$ is a Hermite polynomial of order $q\geq 1$. The first few Hermite polynomials
are $H_1(x)=x$, $H_2(x)=x^2-1$ and $H_3(x)=x^3-3x$.
Why Hermite polynomials? This is because they can be expressed as multiple stochastic integrals.
and because the limit distribution, as $n\to\infty$, of the partial sums $\sum_{k=0}^{n-1}H_q(Y_k)$,
adequably normalized, is known (see Breuer and Major \cite{BM}, Dobrushin and Major \cite{DoMa}, Giraitis and Surgailis \cite{GS} 
and Taqqu \cite{T}). This limit can be Gaussian or not, depending on the order 
$q$ of the polynomial and on the Hurst index
$H$. 

Using our multiple stochastic integral criterion for ASCLT and results from Malliavin calculus, we will show
that the normalized sums of $\sum_{k=0}^{n-1}H_q(Y_k)$, when $q\geq 2$, satisfy an ASCLT if their limit
is Gaussian, that is, when the Hurst index $H$ satisfies $0<H\leq 1-\frac1{2q}$.

Our criterion, when applied to the simple linear case $q=1$, yields the
following ASCLT for fractional Brownian motion. We show that, almost surely,
for any bounded
continuous function $\varphi:\R\to\R$ and any $0<H<1$,
\begin{equation}\label{e:fbmr}
\frac1{\log n}\sum_{k=1}^n \frac1k \varphi\left(\frac{B^H_k}{k^H}\right)
\longrightarrow E[\varphi(N)]\quad\mbox{as $n\to\infty$},
\end{equation}
where $N\sim\NN(0,1)$. Berkes and Horv\'ath \cite{BerkesHorvath},
using other techniques, have obtained a continuous-time version
of (\ref{e:fbmr}).

The paper is organized as follows. In Section 2,
we present the basic elements of Gaussian analysis and Malliavin calculus
used in this paper. The ASCLT criterion is stated and proved in Section 3.
Its application to partial sums of Hermite polynomials of increments of
fractional Brownian motion is presented in Section 4, when the limit
in distribution is Gaussian. In section 5, we discuss the case where the
limit in distribution is non-Gaussian.

\section{Multiple stochastic integrals and Malliavin calculus}\label{sec22}
We shall now present the basic elements of Gaussian analysis and Malliavin calculus that are used in this paper. The reader
is referred to the monograph by  Nualart \cite{Nbook} for any unexplained definition or
 result.

Let $\EuFrak H$ be a real separable Hilbert space. For any $q\geq 1$, let $\EuFrak H^{\otimes q}$ be the $q$th tensor product of $\EuFrak H$ and denote
by $\EuFrak H^{\odot q}$ the associated $q$th symmetric tensor product. We write $X=\{X(h),h\in \EuFrak H\}$ to indicate
an isonormal Gaussian process over
$\EuFrak H$, defined on some probability space $(\Omega ,\mathcal{F},P)$.
This means that $X$ is a centered Gaussian family, whose covariance is given in terms of the
inner product of $\EuFrak H$ by $E\left[ X(h)X(g)\right] =\langle h,g\rangle _{\EuFrak H}$. 

For every $q\geq 1$, let $\mathcal{H}_{q}$ be the $q$th Wiener chaos of $X$,
that is, the closed linear subspace of $L^2(\Omega ,\mathcal{F},P)$
generated by the random variables of the type $\{H_{q}(X(h)),h\in \EuFrak H,\left\|
h\right\| _{\EuFrak H}=1\}$, where $H_{q}$ is the $q$th Hermite polynomial
defined as 
\begin{equation}\label{hermite}
H_q(x) = (-1)^q e^{\frac{x^2}{2}}
 \frac{d^q}{dx^q} \big( e^{-\frac{x^2}{2}} \big).
\end{equation}
We write by convention $\mathcal{H}_{0} = \mathbb{R}$ and $I_{0}(x)=x$, $x\in\mathbb{R}$. For
any $q\geq 1$, the mapping $I_{q}(h^{\otimes q})=q!H_{q}(X(h))$ can be extended to a
linear isometry between \ the symmetric tensor product $\EuFrak H^{\odot q}$
equipped with the modified norm $
\left\| \cdot \right\| _{\EuFrak H^{\odot q}}=
\sqrt{q!}
\left\| \cdot \right\| _{\EuFrak H^{\otimes q}}
$ and the $q$th Wiener chaos $\mathcal{H}_{q}$. 
Then
\[
E[I_p(f)I_q(g)]=\delta_{p,q}\times p!\langle f,g\rangle_{\HH^{\otimes p}}
\quad\mbox{($\delta_{p,q}$ stands for the Kronecker symbol)}
\]
for $f\in\HH^{\odot p}$, $g\in\HH^{\odot q}$ and $p,q\geq 1$.
Moreover, if $f\in\HH^{\otimes q}$, we have
\begin{equation}\label{e:sym}
I_q(f)=I_q(\widetilde{f}),
\end{equation}
where $\widetilde{f}\in\HH^{\odot q}$ is the symmetrization of $f$.

Let $\{e_{k},\,k\geq 1\}$ be a complete orthonormal system in $\EuFrak H$.
Given $f\in \EuFrak H^{\odot p}$ and $g\in \EuFrak H^{\odot q}$, for every
$r=0,\ldots ,p\wedge q$, the \textit{contraction} of $f$ and $g$ of order $r$
is the element of $\EuFrak H^{\otimes (p+q-2r)}$ defined by
\begin{equation}
f\otimes _{r}g=\sum_{i_{1},\ldots ,i_{r}=1}^{\infty }\langle
f,e_{i_{1}}\otimes \ldots \otimes e_{i_{r}}\rangle _{\EuFrak H^{\otimes
r}}\otimes \langle g,e_{i_{1}}\otimes \ldots \otimes e_{i_{r}}
\rangle_{\EuFrak H^{\otimes r}}.  \label{v2}
\end{equation}
Since $f\otimes _{r}g$ is not necessarily symmetric, we denote its
symmetrization by $f\widetilde{\otimes }_{r}g\in \EuFrak H^{\odot (p+q-2r)}$.
Observe that $f\otimes _{0}g=f\otimes g$ equals the tensor product of $f$ and
$g$ while, for $p=q$, $f\otimes _{q}g=\langle f,g\rangle _{\EuFrak H^{\otimes q}}$, namely
the scalar product of $f$ and $g$.
In the particular case $\EuFrak H=L^2(A,\mathcal{A},\mu )$, where
$(A,\mathcal{A})$ is a measurable space and $\mu $ is a $\sigma $-finite and
non-atomic measure, one has that $\EuFrak H^{\odot q}=L_{s}^{2}(A^{q},
\mathcal{A}^{\otimes q},\mu ^{\otimes q})$ is the space of symmetric and
square integrable functions on $A^{q}$. In this case, (\ref{v2}) can be written as
\begin{eqnarray*}
(f\otimes _{r}g)(t_1,\ldots,t_{p+q-2r})
&=&\int_{A^{r}}f(t_{1},\ldots ,t_{p-r},s_{1},\ldots ,s_{r}) \\
&&\times \,g(t_{p-r+1},\ldots ,t_{p+q-2r},s_{1},\ldots ,s_{r})d\mu
(s_{1})\ldots d\mu (s_{r}),
\end{eqnarray*}
that is, we identify $r$ variables in $f$ and $g$ and integrate them out.
The following useful {\sl multiplication formula} holds: if $f\in \EuFrak
H^{\odot p}$ and $g\in \EuFrak
H^{\odot q}$, then
\begin{eqnarray}\label{multiplication}
I_p(f) I_q(g) = \sum_{r=0}^{p \wedge q} r! {p \choose r}{ q \choose r} I_{p+q-2r} (f\widetilde{\otimes}_{r}g).
\end{eqnarray}
\smallskip

Let us now introduce some basic elements of the Malliavin calculus with respect
to the isonormal Gaussian process $X$. Let $\mathcal{S}$
be the set of all
cylindrical random variables of
the form
\begin{equation}
F=\varphi\left( X(h_{1}),\ldots ,X(h_{n})\right) ,  \label{v3}
\end{equation}
where $n\geq 1$, $\varphi:\mathbb{R}^{n}\rightarrow \mathbb{R}$ is an infinitely
differentiable function with compact support and $h_{i}\in \EuFrak H$.
The {\sl Malliavin derivative}  of $F$ with respect to $X$ is the element of
$L^2(\Omega ,\EuFrak H)$ defined as
\begin{equation*}
DF\;=\;\sum_{i=1}^{n}\frac{\partial \varphi}{\partial x_{i}}\left( X(h_{1}),\ldots ,X(h_{n})\right) h_{i}.
\end{equation*}
By iteration, one can
define the $m$th derivative $D^{m}F$, which is an element of $L^2(\Omega ,\EuFrak H^{\odot m})$,
for every $m\geq 2$.
For $m\geq 1$ and $p\geq 1$, ${\mathbb{D}}^{m,p}$ denotes the closure of
$\mathcal{S}$ with respect to the norm $\Vert \cdot \Vert _{m,p}$, defined by
the relation
\begin{equation*}
\Vert F\Vert _{m,p}^{p}\;=\;E\left[ |F|^{p}\right] +\sum_{i=1}^{m}E\left(
\Vert D^{i}F\Vert _{\EuFrak H^{\otimes i}}^{p}\right) .
\end{equation*}
In particular, $DX(h)=h$ for every $h\in \EuFrak H$. The Malliavin derivative $D$ verifies moreover
the following \textsl{chain rule}. If
$\varphi :\mathbb{R}^{n}\rightarrow \mathbb{R}$ is continuously
differentiable with bounded partial derivatives and if $F=(F_{1},\ldots
,F_{n})$ is a vector of elements of ${\mathbb{D}}^{1,2}$, then $\varphi
(F)\in {\mathbb{D}}^{1,2}$ and
\begin{equation*}
D\varphi (F)=\sum_{i=1}^{n}\frac{\partial \varphi }{\partial x_{i}}(F)DF_{i}.
\end{equation*}

Let now $\HH=
\mathrm{L}^{2}(A,\mathcal{A},\mu )$ with $\mu$ non-atomic. Then 
an element $u\in\HH$ can be expressed as $u=\{u_t,\,t\in A\}$ and
the Malliavin
derivative of a multiple integral $F$ of the form $I_q(f)$ (with $f\in\HH^{\odot q}$) 
is
the element
$DF=\{D_tF,\,t\in A\}$
 of $L^2(A\times \Omega )$ given by
\begin{equation}
D_{t}F=D_t \big[I_q(f)\big]=qI_{q-1}\left( f(\cdot ,t)\right).  \label{dtf}
\end{equation}
Thus the derivative of the random variable $I_q(f)$ is the stochastic process $qI_{q-1}\big(f(\cdot,t)\big)$,
$t\in A$.
Moreover, 
\[
\|D\big[I_q(f)\big]\|^2_\HH=q^2\int_A I_{q-1}\left( f(\cdot ,t)\right)^2\mu(dt).
\]

We shall also use the following bound, proved by Nourdin and Peccati in \cite{NP07}, 
for the difference between the law of 
multiple integrals of order $q\geq 2$ with unit variance and the law
of a standard Gaussian random variable.

\begin{prop} \label{noupec}
Let $q\geq 2$ be an integer, $f\in\EuFrak H^{\odot q}$ with
$q!\|f\|^2_{\EuFrak H^{\otimes q}}=E\big[I_q(f)^2\big]=1$
and $N\sim\mathscr{N}(0,1)$.
Then, for all $h:\C\to\R$ such that 
\begin{equation}\label{e:lip}
|h(x)-h(y)|\leq|x-y|,\quad x,y\in\C,
\end{equation}
we have
\begin{equation}\label{was}
\big|E[h(I_q(f))]-E[h(N)]\big|
\leq
\sqrt{E\left[\left(1-\frac1q\|D[I_q(f)]\|^2_\HH\right)^2\right]}.
\end{equation}
\end{prop}

\section{A multiple stochastic integral criterion for ASCLT}\label{gen-crit}

The following result, due to Ibragimov and Lifshits \cite{IL} (Theorem 1.1 therein), 
gives a sufficient condition for extending 
convergence in law to an almost sure limit theorem.

\begin{prop}\label{thm-IL}
Let $\{G_n\}$ be a sequence of random variables converging in distribution towards a random 
variable $G_\infty$, and set
\[
\Delta_n(t)=\frac1{\log n}\sum_{k=1}^n \frac{1}k \big(e^{itG_k}-E(e^{itG_\infty})\big).
\]
If 
\begin{equation}\label{cond-IL}
\sup_{|t|\leq r}\sum_n \frac{E\vert \Delta_n(t)\vert^2}{n\log n}<\infty\quad\mbox{for all $r>0$},
\end{equation}
then, 
almost surely, for all continuous and bounded functions 
$\varphi:\R\to\R$, we have
\[
\frac{1}{\log n}\sum_{k=1}^{n} \frac{1}{k}\,\varphi(G_k) 
\longrightarrow E[\varphi(G_\infty)]\quad\mbox{as $n\to\infty$}. 
\]
\end{prop}

The following theorem provides a criterion for an ASCLT for multiple stochastic integrals. 
It is expressed in terms of the kernels of these integrals. 

\begin{thm}\label{main}
Let the notation of Section \ref{sec22} prevail.
Fix $q\geq 2$, and let $\{G_n\}$ be a sequence of the form $G_n=I_q(f_n)$, with $f_n\in\EuFrak H^{\odot q}$. Assume that
$E[G_n^2]=q!\|f_n\|^2_{\HH^{\otimes q}}=1$ for all $n$, and that $G_n\overset{\rm law}{\longrightarrow}N\sim
\NN(0,1)$ as $n\to\infty$.
If the following two conditions are satisfied:
\begin{eqnarray*}
(A_1) &\quad&\displaystyle{\sum_{n\geq 2} \frac1{n\log^2n}\sum_{k=1}^n \frac1k\,
\|f_k\otimes_r f_k\|_{\EuFrak H^{\otimes 2(q-r)}}}<\infty\quad\mbox{for every $r=1,\ldots,q-1$};\\
(A_2) &\quad&\displaystyle{\sum_{n\geq 2} \frac1{n\log^3n}\sum_{k,l=1}^n \frac{\big|
\langle f_k,f_l\rangle_{\EuFrak H^{\otimes q}}
\big|}
{kl}}<\infty,
\end{eqnarray*}
then $\{G_n\}$ satisfies an ASCLT. In other words, almost surely, for all 
continuous and bounded $\varphi:\R\to\R$, 
\[
\frac{1}{\log n}\sum_{k=1}^{n} \frac{1}{k}\,\varphi(G_k) \,\longrightarrow\, E[\varphi(N)],\quad\mbox{as $n\to\infty$}.
\]
\end{thm}
\begin{rem}
{\rm
We have $E\big(G_kG_l\big)=q!\langle f_k,f_l\rangle_{\EuFrak H^{\otimes q}}$. Consequently, condition $(A_2)$
can be replaced by the following equivalent condition:
\[
(A_2') \quad\quad\displaystyle{\sum_{n\geq 2} \frac1{n\log^3n}\sum_{k,l=1}^n \frac{\big|
E(G_kG_l)
\big|}
{kl}}<\infty.
\]
}
\end{rem}
{\bf Proof of Theorem \ref{main}}. We shall verify the sufficient condition (\ref{cond-IL}),
that is the Ibragimov-Lifshits criterion. For simplicity, let $g(t)=E(e^{itN})=e^{-t^2/2}$.
Then
\begin{eqnarray*}
&&E\vert \Delta_n(t)\vert^2 \\&=&\frac{1}{\log^2 n}\sum_{k,l=1}^n \frac{1}{kl}
E\left[\big(e^{itG_k}-g(t)\big)\big(e^{-itG_l}-g(t)\big)\right]\\
&=&\frac{1}{\log^2 n} \sum_{k,l=1}^{n}\frac{1}{kl}
\left[
E\big(e^{it(G_k-G_l)}\big)-g(t)\left(E\big(e^{itG_k}\big)+E\big(e^{-itG_l}\big)\right)+g^2(t)
\right]\\
&=&\frac{1}{\log^2 n} \sum_{k,l=1}^{n}\frac{1}{kl}
\left[
\left(E\big(e^{it(G_k-G_l)}\big)-g^2(t)\right) 
- g(t)\left(
E\big(e^{itG_k}\big)-g(t)
\right)
- g(t)\left(
E\big(e^{-itG_l}\big)-g(t)
\right)
\right].
\end{eqnarray*}
Let $t\in\R$ and $r>0$ be such that $|t|\leq r$.
The function $\varphi(x)=\frac1r e^{itx}$
satisfies $|\varphi'(x)|\leq 1$ and hence (\ref{e:lip}). Therefore, 
by (\ref{was}), we have, for $j\in\{k,l\}$,
\[
\left| E\big(e^{\pm itG_j}\big) - g(t)\right| \leq r
\sqrt{E\left[\left(1-\frac1q\|DG_j\|^2_\HH\right)^2\right]}.
\]
Similarly,
\begin{eqnarray*}
\left| E\big(e^{it(G_k-G_l)}\big) - g^2(t)\right|&=&\left| E\big(e^{it\sqrt{2}\,\frac{G_k-G_l}{\sqrt{2}}}\big) - g(\sqrt{2}\,t)\right| \notag\\
&\leq& \sqrt{2}\,r \sqrt{
E\left[\left(1-\frac1q\left\|D\big(\frac{G_k-G_l}{\sqrt{2}}\big) \right\|^2_\HH\right)^2\right]}\\
&=& \sqrt{2}\,r \sqrt{
E\left[\left(1-\frac1{2q}\|DG_k-DG_l\|^2_\HH\right)^2\right]}.
\end{eqnarray*}
But
\begin{eqnarray*}
1-\frac1{2q}\|DG_k-DG_l\|^2_\HH
&=&\frac12\left(
1-\frac1q\|DG_k\|^2_\HH
\right)+\frac12\left(
1-\frac1q\|DG_l\|^2_\HH
\right)+\frac1q\langle DG_k,DG_l\rangle_\HH
\end{eqnarray*}
so that, since $(x+y+z)^2\leq 3(x^2+y^2+z^2)$ and $\sqrt{u+v+w}\leq \sqrt{u}+\sqrt{v}+\sqrt{w}$,
\begin{eqnarray*}
\sqrt{E\left[\left(1-\frac1q\|DG_k-DG_l\|^2_\HH\right)^2\right]}
&\leq&\frac{\sqrt{3}}2\sqrt{E\left[\left(
1-\frac1q\|DG_k\|^2_\HH
\right)^2\right]}
\\
&&+\frac{\sqrt 3}{2}\sqrt{E\left[\left(
1-\frac1q\|DG_l\|^2_\HH
\right)^2\right]}
+\frac{\sqrt{3}}q\sqrt{E\big[\langle DG_k,DG_l\rangle_\HH^2\big]}.
\end{eqnarray*}
Consequently, to get (\ref{cond-IL}), it suffices to prove
that the two following conditions hold:
\begin{equation}
\sum_{n}\frac{1}{n\log^3n}\sum_{k,l=1}^{n}\frac{1}{kl}
\sqrt{E\left(
1-\frac1q\|DG_k\|^2_\HH
\right)^2}<\infty\label{cond1bis}
\end{equation}
and
\begin{equation}
\sum_{n}\frac{1}{n\log^3n}\sum_{k,l=1}^{n}\frac{1}{kl}
\sqrt{E\langle DG_k,DG_l\rangle_\HH^2}<\infty\label{cond2}.
\end{equation}
Since $\sum_{l=1}^n \frac1l \leq 1+\log n$, one can observe that (\ref{cond1bis}) is 
a consequence of
\begin{equation}
\sum_{n}\frac{1}{n\log^2n}\sum_{k=1}^{n}\frac{1}{k}
\sqrt{E\left(
1-\frac1q\|DG_k\|^2_\HH
\right)^2}
<\infty.\label{cond1}
\end{equation}
We shall now prove that (\ref{cond1}) and (\ref{cond2}) are satisfied.\\

{\it Proof of (\ref{cond1})}. In view of 
the normalization $q!\|f_k\|^2_\HH=1$ for all $k$, we have, by Lemma \ref{key-lm} below:
\[
\frac1q\|DG_k\|^2_{\HH}-1 = q\sum_{r=1}^{q-1} (r-1)!\binom{q-1}{r-1}^2 I_{2q-2r}(f_k\widetilde{\otimes}_r f_k).
\]
Hence, taking into account the orthogonality between multiple stochastic integrals, we have
\[
E\left[\left(1-\frac1q\|DG_k\|^2_{\HH} \right)^2\right] = q^2\sum_{r=1}^{q-1} (r-1)!^2\binom{q-1}{r-1}^4 
(2q-2r)! \,\|f_k\widetilde{\otimes}_r f_k\|^2_{\HH^{\otimes(2q-2r)}}.
\]
Via the straightforward inequality 
\begin{equation}\label{ineqsquaroot}
\sqrt{x_1^2+\ldots+x_{q}^2}\leq |x_1|+\ldots+|x_{q}|,
\end{equation} 
and since $\|f_k\widetilde{\otimes}_rf_k\|_{\HH^{\otimes(2q-2r)}}  \leq 
\|f_k\otimes_rf_k\|_{\HH^{\otimes(2q-2r)}}$, we get
\[
\sqrt{E\left[\left(1-\frac1q\|DG_k\|^2_{\HH} \right)^2\right]} \leq 
q\sum_{r=1}^{q-1} (r-1)!\binom{q-1}{r-1}^2 \sqrt{(2q-2r)!} 
\,\|f_k\otimes_r f_k\|_{\HH^{\otimes(2q-2r)}}.
\]
Combining all these bounds, we obtain
\begin{eqnarray*}
&&\sum_{n}\frac{1}{n\log^2 n}\sum_{k=1}^n \frac1k 
\sqrt{E\left(
1-\frac1q\|DG_k\|^2_\HH
\right)^2}\\
&\leq&q\sum_{r=1}^{q-1}
(r-1)!\binom{q-1}{r-1}^2 \sqrt{(2q-2r)!} 
\times
\sum_{n}\frac{1}{n\log^2 n}\sum_{k=1}^n \frac1k 
\,\|f_k\otimes_r f_k\|_{\HH^{\otimes(2q-2r)}},
\end{eqnarray*}
so that assumption
$(A_1)$ immediately implies (\ref{cond1}).\\

{\it Proof of (\ref{cond2})}. By Lemma \ref{key-lm} below and the orthogonality between multiple stochastic integrals, 
we have
\[
E\langle DG_k,DG_l\rangle^2_{\HH} = q^4 \sum_{r=1}^{q-1} (r-1)!^2\binom{q-1}{r-1}^4 (2q-2r)!\,
\|f_k\otimes_r f_l\|^2
_{\HH^{\otimes(2q-2r)}}
+(q\,q!)^2\big| \langle f_k, f_l\rangle
_{\HH^{\otimes q}}\big|^2.
\]
By Lemma \ref{contr} below, we also have
\[
\|f_k\otimes_r f_l\|^2
_{\HH^{\otimes(2q-2r)}}=
\langle f_k\otimes_{q-r}f_k,f_l\otimes_{q-r} f_l\rangle
_{\HH^{\otimes 2r}},
\]
so that, by Cauchy-Schwarz inequality:
\begin{eqnarray*}
\|f_k\otimes_r f_l\|^2
_{\HH^{\otimes(2q-2r)}}&\leq&
\|f_k\otimes_{q-r} f_k\|
_{\HH^{\otimes 2r}}
\|f_l\otimes_{q-r} f_l\|
_{\HH^{\otimes 2r}}\\
&\leq&\frac12\left(
\|f_k\otimes_{q-r} f_k\|
_{\HH^{\otimes 2r}}^2
+
\|f_l\otimes_{q-r} f_l\|
_{\HH^{\otimes 2r}}^2
\right).
\end{eqnarray*}
Consequently, using again (\ref{ineqsquaroot}), we obtain that
\begin{eqnarray*}
\sqrt{E\langle DG_k,DG_l\rangle^2_{\HH}} &\leq& 
q^2 \sum_{r=1}^{q-1} (r-1)!\binom{q-1}{r-1}^2 \sqrt{(2q-2r)!}\\
&&\times
\frac{1}{\sqrt{2}}\left(
\|f_k\otimes_{q-r} f_k\|
_{\HH^{\otimes 2r}}
+
\|f_l\otimes_{q-r} f_l\|
_{\HH^{\otimes 2r}}
\right)+q\,q!\big| \langle f_k, f_l\rangle
_{\HH^{\otimes q}}\big|.
\end{eqnarray*}
Finally, we obtain (\ref{cond2}) from the conjunction of $(A_1)$ and $(A_2)$, which
completes the proof of Theorem \ref{main}.
\fin
In the previous proof, we used the two following lemmas:
\begin{lemma}\label{key-lm}
Consider two random variables $F=I_q(f)$, $G=I_q(g)$, with $f,g\in\HH^{\odot q}$. Then
\[
\langle DF,DG\rangle_{\HH} = q^2 \sum_{r=1}^{q-1} (r-1)!\binom{q-1}{r-1}^2 I_{2q-2r}(f\widetilde{\otimes}_r g) 
+ q\,q!\langle f,g\rangle_{\HH^{\otimes q}}.
\]
\end{lemma}
{\bf Proof}.
Without loss of generality, we can assume that $\HH=L^{2}(A,\mathscr{A},\mu)$ where
$(A,\mathscr{A})$ is a measurable space, and $\mu$ is a $\sigma$-finite and non-atomic measure.
Thus, we can write
\begin{eqnarray*}
\langle DF,DG\rangle_\HH &=&\int_A D_tF\,D_tG\,\mu(dt)\\
&=&q^2\int_A I_{q-1}\big(f(\cdot,t)\big)I_{q-1}\big(g(\cdot,t)\big)\mu(dt)\quad\mbox{by (\ref{dtf})}\\
&=&q^2\int_A \sum_{r=0}^{q-1} r!\binom{q-1}{r}^2 I_{2q-2-2r}\big(f(\cdot,t)\widetilde{\otimes}_r g(\cdot,t)\big)\mu(dt)\quad\mbox{by (\ref{multiplication})}\\
&=&q^2\int_A \sum_{r=0}^{q-1} r!\binom{q-1}{r}^2 I_{2q-2-2r}\big(f(\cdot,t)\otimes_r g(\cdot,t)\big)\mu(dt)\quad\mbox{by (\ref{e:sym})}
\end{eqnarray*}
\begin{eqnarray*}
&=&q^2\sum_{r=0}^{q-1} r!\binom{q-1}{r}^2 I_{2q-2-2r}(f\otimes_{r+1}g)\quad\mbox{by linearity}\\
&=&q^2\sum_{r=0}^{q-1} r!\binom{q-1}{r}^2 I_{2q-2-2r}(f\widetilde{\otimes}_{r+1}g)\\
&=&q^2 \sum_{r=1}^{q} (r-1)!\binom{q-1}{r-1}^2 I_{2q-2r}(f\widetilde{\otimes}_r g)\\
&=&q^2 \sum_{r=1}^{q-1} (r-1)!\binom{q-1}{r-1}^2 I_{2q-2r}(f\widetilde{\otimes}_r g) + q\,q!\langle f,g\rangle_{\HH^{\otimes q}}.
\end{eqnarray*}
\fin

\begin{lemma}\label{contr}
Let $f,g\in\HH^{\odot q}$. Then
$\|f\otimes_r g\|^2_{\HH^{\otimes (2q-2r)}}=
\langle f\otimes_{q-r}f,g\otimes_{q-r}g\rangle_{\HH^{\otimes 2r}}
.$
\end{lemma}
{\bf Proof}.
Without loss of generality, we can assume that $\HH=L^{2}(A,\mathscr{A},\mu)$ where
$(A,\mathscr{A})$ is a measurable space, and $\mu$ is a $\sigma$-finite and non-atomic measure.
Using the definition of contractions and Fubini theorem, we can write
\begin{eqnarray*}
&&\|f\otimes_r g\|^2_{\HH^{\otimes (2q-2r)}}\\
&=&
\int_{A^{2q-2r}}\mu(dx_1)\ldots \mu(dx_{q-r})\mu(dy_1)\ldots \mu(dy_{q-r})\\
&&\times\left(
\int_{A^r}\mu(dz_1)\ldots \mu(dz_r)\,f(x_1,\ldots,x_{q-r},z_1,\ldots,z_r)g(y_1,\ldots,y_{q-r},z_1,\ldots,z_r)
\right)^2\\
&=&
\int_{A^{2q-2r}}\mu(dx_1)\ldots \mu(dx_{q-r})\mu(dy_1)\ldots \mu(dy_{q-r})\\
&&\times
\int_{A^r}\mu(dz_1)\ldots \mu(dz_r)\,f(x_1,\ldots,x_{q-r},z_1,\ldots,z_r)g(y_1,\ldots,y_{q-r},z_1,\ldots,z_r)\\
&&\times
\int_{A^r}\mu(dt_1)\ldots \mu(dt_r)\,f(x_1,\ldots,x_{q-r},t_1,\ldots,t_r)g(y_1,\ldots,y_{q-r},t_1,\ldots,t_r)\\
&=&
\int_{A^{2r}}\mu(dz_1)\ldots \mu(dz_r)
\mu(dt_1)\ldots \mu(dt_r)\\
&&\times
\int_{A^{q-r}}
\mu(dx_1)\ldots \mu(dx_{q-r})
\,f(x_1,\ldots,x_{q-r},z_1,\ldots,z_r)f(x_1,\ldots,x_{q-r},t_1,\ldots,t_r)\\
&&\times
\int_{A^{q-r}}
\mu(dy_1)\ldots \mu(dy_{q-r})
g(y_1,\ldots,y_{q-r},z_1,\ldots,z_r)\,g(y_1,\ldots,y_{q-r},t_1,\ldots,t_r)\\
&=&
\int_{A^{2r}}\mu(dz_1)\ldots \mu(dz_r)
\mu(dt_1)\ldots \mu(dt_r)\\
&&\times
f\otimes_{q-r}f(z_1,\ldots,z_r,t_1,\ldots,t_r)
g\otimes_{q-r}g(z_1,\ldots,z_r,t_1,\ldots,t_r)
=\langle f\otimes_{q-r}f,g\otimes_{q-r}g\rangle_{\HH^{\otimes 2r}}.
\end{eqnarray*}
\fin

When $q=1$, we have the following result:
\begin{prop}\label{q=1}
Let $\{G_n\}$ be a centered Gaussian sequence with unit variance. If the following condition is satisfied,
\[
(A_2') \quad\quad\displaystyle{\sum_{n\geq 2} \frac1{n\log^3n}\sum_{k,l=1}^n \frac{\big|
E(G_kG_l)
\big|}
{kl}}<\infty,
\]
then $\{G_n\}$ satisfies an ASCLT. In other words, almost surely, for all 
continuous and bounded $\varphi:\R\to\R$,
\[
\frac{1}{\log n}\sum_{k=1}^{n} \frac{1}{k}\,\varphi(G_k) \,
\longrightarrow\, E[\varphi(N)],\quad\mbox{as $n\to\infty$}. 
\]
\end{prop}
{\bf Proof}. Let $t\in\R$ and $r>0$ be such that $|t|\leq r$.
Since $E[e^{itG_k}]$ equals $g(t)=e^{-t^2/2}$, we have
\begin{eqnarray*}
E\vert \Delta_n(t)\vert^2 &=&\frac{1}{\log^2 n}\sum_{k,l=1}^n \frac{1}{kl}
E\left[\big(e^{itG_k}-e^{-t^2/2}\big)\big(e^{-itG_l}-e^{-t^2/2}\big)\right]\\
&=&\frac{1}{\log^2 n} \sum_{k,l=1}^{n}\frac{1}{kl}
\left[
E\big(e^{it(G_k-G_l)}\big)-e^{-t^2}
\right]\\
&=&\frac{1}{\log^2 n} \sum_{k,l=1}^{n}\frac{e^{-t^2}}{kl}
\big(e^{E(G_kG_l)t^2}-1\big)\\
&\leq&\frac{r^2e^{r^2}}{\log^2 n} \sum_{k,l=1}^{n}\frac{\big|E(G_kG_l)\big|}{kl}\quad\mbox{since
$|e^{x}-1|\leq e^{|x|}|x|$ and $|E(G_kG_l)|\leq 1$}.
\end{eqnarray*}
Therefore, assumption $(A'_2)$ implies (\ref{cond-IL}), and the proof of the proposition is done. 
\fin

\section{Partial sums of Hermite polynomials of increments of fBm: the Gaussian case}

Let $B^H$ be a
 fractional Brownian motion 
(fBm) with Hurst index $H\in(0,1)$. 
We are interested in an ASCLT for the $q$-Hermite power variations of $B^H$,  defined as
\begin{equation}\label{her-intro}
V_n=\sum_{k=0}^{n-1} H_q(B^H_{k+1}-B^H_{k}),\quad n\geq 1,
\end{equation}
in cases where $V_n$, adequably normalized, converges to a normal distribution.
Here, $H_q$ stands for the Hermite polynomial of degree $q$, given by 
(\ref{hermite}).

We first treat the case $q=1$ and $0<H<1$.
Convergence in distribution of
\[
G_n=\frac{V_n}{n^H}=\frac{B^H_n}{n^H}
\]
to a normal law is trivial because, by self-similarity,
$G_n\overset{\rm law}{=}B^H_1\sim\NN(0,1)$.
The following theorem provides the corresponding ASCLT. A continuous
time version of the result was obtained by 
Berkes and Horv\'ath \cite{BerkesHorvath}
using different methods.

\begin{thm}\label{thm1234}
For all $H\in(0,1)$, the sequence $\{G_n\}$ satisfies an ASCLT. In other words,
almost surely, for all continuous and bounded $\varphi:\R\to\R$,
\[
\frac1{\log n}\sum_{k=0}^{n-1}\frac1k\,\varphi(B^H_k/k^H)\longrightarrow E[\varphi(N)]\quad\mbox{as $n\to\infty$}.
\]
\end{thm}
{\bf Proof}.
We shall apply Proposition \ref{q=1}. The cases $H<1/2$ and $H\geq 1/2$ are treated separately.
From now on, the value of a constant $C>0$ may change from line to
line.\\

{\it Case $H<1/2$.} For any $b\geq a\geq 0$, we have
\[
b^{2H}-a^{2H}=2H \int_0^{b-a}\frac{dx}{(x+a)^{1-2H}}
\leq 2H\int_0^{b-a}\frac{dx}{x^{1-2H}}=(b-a)^{2H}.
\]
Hence, for $l\geq k\geq 1$, we have $l^{2H}-(l-k)^{2H}\leq k^{2H}$ so that
\[
|E[B^H_kB^H_l]|=\frac12\big(k^{2H}+l^{2H}-(l-k)^{2H}\big)\leq k^{2H}.
\]
Thus
\begin{eqnarray*}
\sum_{n\geq 2}\frac1{n\log^3 n}\sum_{l=1}^n\frac1l\sum_{k=1}^l\frac{|E[G_kG_l]|}{k}
&=&\sum_{n\geq 2}\frac1{n\log^3 n}\sum_{l=1}^n\frac1{l^{1+H}}\sum_{k=1}^l\frac{|E[B^H_kB^H_l]|}{k^{1+H}}\\
&\leq&\sum_{n\geq 2}\frac1{n\log^3 n}\sum_{l=1}^n\frac1{l^{1+H}}\sum_{k=1}^l\frac1{k^{1-H}}\\
&\leq&C\sum_{n\geq 2}\frac1{n\log^3 n}\sum_{l=1}^n\frac1{l}
\leq C\sum_{n\geq 2}\frac1{n\log^2 n}<\infty.
\end{eqnarray*}
Consequently, condition $(A'_2)$ in Proposition \ref{q=1} is fulfilled.\\

{\it Case $H\geq 1/2$}. For $l\geq k\geq 1$, it follows from (\ref{e:cov})
that
\begin{eqnarray*}
|E[B^H_kB^H_l]|&=&\left| 
\sum_{i=0}^{k-1}\sum_{j=0}^{l-1}
E[(B^H_{i+1}-B^H_i)(B^H_{j+1}-B^H_j)]
\right|
\leq\sum_{i=0}^{k-1}\sum_{j=0}^{l-1}|\rho(i-j)|\\
&\leq&k\sum_{r=-l+1}^{l-1}|\rho(r)|\leq Ckl^{2H-1}.
\end{eqnarray*}
The last inequality comes from the fact that 
$\rho(0)=1$, $\rho(1)=\rho(-1)=(2^{2H}-1)/2$
and, if $r\geq 2$,
\begin{eqnarray*}
|\rho(-r)|&=&|\rho(r)|=\big| E[(B^H_{r+1}-B^H_r)B^H_1]
=H(2H-1)\int_0^1 du\int_r^{r+1}dv(v-u)^{2H-2}\\
&\leq& H(2H-1)\int_0^1 (r-u)^{2H-2}du \leq H(2H-1)(r-1)^{2H-2}.
\end{eqnarray*} 
Consequently, 
\begin{eqnarray*}
\sum_{n\geq 2}\frac1{n\log^3 n}\sum_{l=1}^n\frac1l\sum_{k=1}^l\frac{|E[G_kG_l]|}{k}
&=&\sum_{n\geq 2}\frac1{n\log^3 n}\sum_{l=1}^n\frac1{l^{1+H}}\sum_{k=1}^l\frac{|E[B^H_kB^H_l]|}{k^{1+H}}\\
&\leq&C\sum_{n\geq 2}\frac1{n\log^3 n}\sum_{l=1}^n\frac1{l^{2-H}}\sum_{k=1}^l\frac1{k^H}\\
&\leq&C\sum_{n\geq 2}\frac1{n\log^3 n}\sum_{l=1}^n\frac1{l}
\leq C\sum_{n\geq 2}\frac1{n\log^2 n}<\infty.
\end{eqnarray*}
Finally, condition $(A'_2)$ in Proposition \ref{q=1} is satisfied,
which completes the proof of Theorem \ref{thm1234}.
\fin

In the remaining part of this section, we assume that $q\geq 2$ 
and $0< H\leq 1-1/(2q)$. 
When $H\neq 1/2$, since the increments of $B^H$ are not independent 
and $V_n$ is not linear,
the asymptotic behavior of \eqref{her-intro} is more difficult to investigate. 
In fact, thanks to the seminal works of
Breuer and Major \cite{BM}, 
%Dobrushin and Major \cite{DoMa}, 
Giraitis and Surgailis \cite{GS} 
and Taqqu \cite{T}, it is  known that,
as $n\to\infty$:
\begin{itemize}
\item If $q\geq 2$ and $0<H<1-1/(2q)$ then
\begin{equation}\label{eq:Breuer_Major1}
G_n:=\frac{V_n}{ \sigma_n\,\sqrt n}
\,\overset{{\rm Law}}{\longrightarrow}\,
\NN(0,1).
\end{equation}
\item If $q\geq 2$ and $H=1-1/(2q)$ then
\begin{equation}
\label{eq:Breuer_Major2}
G_n:=\frac{V_n}{\sigma_n\sqrt{ n \log n}}
\,\overset{{\rm Law}}{\longrightarrow}\,
\NN(0,1).
\end{equation}
\end{itemize}
Here, $\sigma_n$ denotes the positive normalizing constant
which ensures that $E[G_n^2]=1$. 
The case $H>1-\frac1{2q}$ will be considered in Section \ref{s:hermite}.
Proofs of (\ref{eq:Breuer_Major1}) and (\ref{eq:Breuer_Major2}), together with rates of convergence,
can be found in \cite{NP07} and \cite{BN}, respectively.
 
We want to see if one can associate almost sure central limit theorems to the convergences 
(\ref{eq:Breuer_Major1}) and (\ref{eq:Breuer_Major2}). To do so, we need a few lemmas.
\eject
\begin{lemma}\label{sigma}
Let $q\geq 2$. As $n\to\infty$, 
\begin{enumerate}
\item if $H<1-\frac1{2q}$, then $\sigma_n\to q!\sum_{r\in\mathbb{Z}} \rho(r)^q>0$;
\item if $H=1-\frac1{2q}$, then $\sigma_n\to 2q!\big(1-\frac1{2q}\big)^q\big(1-\frac1q\big)^q>0$.
\end{enumerate}
\end{lemma}
{\bf Proof}.
We have $E[(B^H_{k+1}-B^H_k)(B^H_{l+1}-B^H_l)]=\rho(k-l)$
where $\rho$ is given in (\ref{e:cov}).
Recall that $\rho$ is an even function and that
\[
\rho(r)=H(2H-1)r^{2H-2} + o(r^{2H-2}),\quad\mbox{as $|r|\to\infty$}.
\]
We deduce that
$\sum_{r\in\mathbb{Z}} |\rho(r)|^q<\infty$ if and only if $H<1-1/(2q)$.
On the other hand, 
\begin{eqnarray*}
E[V_n^2]&=&\sum_{k,l=0}^{n-1}E\big(H_q(B^H_{k+1}-B^H_k)H_q(B^H_{l+1}-B^H_l)\big)
=q!\sum_{k,l=0}^{n-1}\rho(k-l)^q\\
&=&q!\sum_{l=0}^{n-1}\sum_{r=-l}^{n-1-l}\rho(r)^q
=q!\sum_{|r|<n} \big(n-1-|r|\big)\rho(r)^q\\
&=&q!\left(n\sum_{|r|<n} \rho(r)^q - \sum_{|r|<n} \big(|r|+1\big)
\rho(r)^q\right).
\end{eqnarray*}
Assume first that $H<1-1/(2q)$. In this case, 
\[
\sigma_n^2=q!\sum_{r\in\mathbb{Z}} \rho(r)^q\left(1-\frac{|r|+1}{n}\right)
{\bf 1}_{\{|r|<n\}}.
\]
In addition, we also have 
$
\sum_{r\in\mathbb{Z}} \vert\rho(r)\vert^q<\infty. 
$
Hence, we deduce by bounded Lebesgue convergence that
\[
\sigma_n^2\to q!\sum_{r\in\mathbb{Z}} \rho(r)^q,\quad\mbox{as $n\to\infty$}.
\]
Assume now that $H=1-\frac1{2q}$. In that case,  as $|r|\to\infty$,
\[
\rho(r)^q\sim
H^q(2H-1)^q |r|^{(2H-2)q}=
\left(1-\frac1{2q}\right)^q\left(1-\frac1q\right)^q\frac{1}{|r|}.
\]
Therefore, as $n\to\infty$,
\[
\sum_{|r|<n} \rho(r)^q \sim
\left(1-\frac1{2q}\right)^q\left(1-\frac1q\right)^q\sum_{0<|r|<n}\frac{1}{|r|}
\sim 2\left(1-\frac1{2q}\right)^q\left(1-\frac1q\right)^q \log n
\]
and 
\[
\sum_{|r|<n} \big(|r|+1\big)\rho(r)^q \sim 
\left(1-\frac1{2q}\right)^q\left(1-\frac1q\right)^q\sum_{|r|<n} 1
\sim 2n
\left(1-\frac1{2q}\right)^q\left(1-\frac1q\right)^q.
\]
The desired conclusion follows.\fin

The next lemma follows from Nourdin and Peccati \cite{NP07} and Breton and Nourdin \cite{BN}. It will play a
crucial role in the proof of Theorem \ref{BM}.
\begin{lemma}\label{Stein-lemma}
Fix $q\geq 2$, and let $\HH$ be the real and separable Hilbert space
defined as follows: (i) denote by $\mathscr{E}$ the
set of all $\mathbb{R}$-valued step functions on $[0,\infty)$, (ii)
define $\EuFrak H$ as the Hilbert space obtained by closing
$\mathscr{E}$ with respect to the scalar product
\[
\left\langle
{\mathbf{1}}_{[0,t]},{\mathbf{1}}_{[0,s]}\right\rangle _{\EuFrak
H}=E[B^H_tB^H_s].
\]
For any $n\geq 2$, let $f_n$ be the element of $\HH^{\odot q}$ defined by
\begin{equation}\label{fk}
f_n=\left\{
\begin{array}{lll}
\frac1{\sigma_n\sqrt{n}}
\sum_{k=0}^{n-1}{\bf 1}_{[k,k+1]}^{\otimes q}&\quad&\mbox{if $H<1-\frac1{2q}$}\\
\\
\frac{1}{\sigma_n\,\sqrt{n\log n}}\sum_{k=0}^{n-1}{\bf 1}_{[k,k+1]}^{\otimes q}&\quad&\mbox{if $H=1-\frac1{2q}$}
\end{array}
\right.
\end{equation}
with $\sigma_n$ the positive normalizing constant which ensures that $G_n=I_q(f_n)$ has variance one. 
Then there exists a constant $C>0$, depending only on $q$ and $H$ (but not on $n$), such that,
for all $n\geq 1$ and $r=1,\ldots,q-1$:
\[
\|f_n\otimes_r f_n\|_{\HH^{\otimes (2q-2r)}}\leq C\times
\left\{
\begin{array}{lll}
n^{-1/2}&\quad&\mbox{if $H\leq \frac12$}\\
\\
n^{H-1}&\quad&\mbox{if $\frac12\leq H\leq \frac{2q-3}{2q-2}$}\\
\\
n^{qH-q+1/2}&\quad&\mbox{if $\frac{2q-3}{2q-2}\leq H< 1-\frac1{2q}$}\\
\\
(\log n)^{-1/2}&\quad&\mbox{if $H= 1-\frac1{2q}$}
\end{array}
\right..
\]
\end{lemma}

We can now state and prove the main result of this section.
\begin{thm}\label{BM}
Let $q\geq 2$ and $H\leq 1-1/(2q)$.
For $n\geq 1$, set
\[
V_n=\sum_{k=0}^{n-1} H_q(B^H_{k+1}-B^H_{k}),
\] 
and define
\[
G_n=\left\{
\begin{array}{lll}
V_n/(\sigma_n\,\sqrt n)&\quad&\mbox{if $H<1-\frac1{2q}$}\\
\\
V_n/(\sigma_n\sqrt{ n \log n})&\quad&\mbox{if $H=1-\frac1{2q}$}
\end{array}
\right..
\] 
Here, $\sigma_n$ denotes the positive normalizing constant
which ensures that $E[G_n^2]=1$.
Then $\{G_n\}$ satisfies an ASCLT. In other words, almost surely, for all continuous and bounded $\varphi:\R\to\R$,
\[
\frac{1}{\log n}\sum_{k=1}^{n} \frac{1}{k}\,\varphi(G_k) 
\longrightarrow E[\varphi(N)]\quad\mbox{as $n\to\infty$}. 
\]
\end{thm}
{\bf Proof}. We shall apply Theorem \ref{main}, and let $C$ be a positive
constant, depending only on $q$ and $H$, whose value changes from line to line. 
We consider 
the real and separable Hilbert space
$\EuFrak H$ as defined in Lemma \ref{Stein-lemma}.

First, we focus on the case $H<1-1/(2q)$. 
We have $G_n=I_q(f_n)$ with $f_n$ given by (\ref{fk}). 
Let us verify assumptions $(A_1)$ and $(A_2)$ in Theorem \ref{main}. 
According to Lemma \ref{Stein-lemma}, 
there exists $\alpha>0$, depending only on $q$ and $H$
(but not on $k$ and $r$), such that, for all $k\geq 1$ and $r=1,\ldots,q-1$,
$\|f_k\otimes_r f_k\|_{\HH^{\otimes (2q-2r)}}\leq Ck^{-\alpha}$.
Hence
\[
\sum_{n\geq 2}\frac{1}{n\log^2 n}\sum_{k=1}^n \frac{1}{k}\|f_k\otimes_r f_k\|_{\HH^{\otimes (2q-2r)}}
\leq
C\sum_{k\geq 1}\frac{1}{k^{1+\alpha}}\times \sum_{n\geq 2}\frac{1}{n\log^2 n}
<\infty,
\]
that is assumption $(A_1)$ is verified.
On the other hand, we have
\[
\langle f_k,f_l\rangle_{\HH^{\otimes q}} =\frac{1}{\sigma_k\sigma_l\,\sqrt{kl}}
\sum_{i=0}^{k-1}\sum_{j=0}^{l-1}\rho(j-i)^q
\]
with $\rho$ given by (\ref{e:cov}).
Since $\sigma_k\rightarrow\sigma_\infty>0$ as $k\to\infty$ (see Lemma \ref{sigma}), we have,
for $l\geq k\geq 1$,
\begin{eqnarray*}
\big|\langle f_k,f_l\rangle_{\HH^{\otimes q}}\big|&\leq&\frac{C}{\sqrt{kl}}
\sum_{i=0}^{k-1}\sum_{j=0}^{l-1}\big|\rho(j-i)\big|^q
=\frac{C}{\sqrt{kl}}\sum_{i=0}^{k-1}\sum_{r=-i}^{l-1-i}\big|\rho(r)\big|^q\\
&\leq&C\sqrt{\frac{k}l}\,\sum_{r\in\Z}\big|\rho(r)\big|^q \leq C\sqrt{\frac{k}{l}},
\end{eqnarray*}
where the last inequality follows from the fact that $\sum_{r\in\Z}\big|\rho(r)\big|^q<\infty$.
Consequently, assumption $(A_2)$ is verified as well, because
\begin{eqnarray*}
\sum_{n\geq 2}\frac{1}{n\log^3 n}\sum_{k,l=1}^{n} \frac{\big|\langle f_k,f_l\rangle_{\HH^{\otimes q}}\big|}{kl}
&\leq&
2\sum_{n\geq 2}\frac{1}{n\log^3 n}\sum_{l=1}^{n}\sum_{k=1}^l \frac{\big|\langle f_k,f_l\rangle_{\HH^{\otimes q}}\big|}{kl}\\
&\leq&C\sum_{n\geq 2}\frac{1}{n\log^3 n}\sum_{l=1}^n \frac{1}{l^{3/2}}\sum_{k=1}^l \frac1{\sqrt{k}}\\
&\leq&C\sum_{n\geq 2}\frac{1}{n\log^3 n}\sum_{l=1}^n\frac{1}{l}
\leq C\sum_{n\geq 2}\frac{1}{n\log^2 n}<\infty.
\end{eqnarray*}

It remains now to consider the critical case $H=1-1/(2q)$. 
We have $G_n=I_q(f_n)$ with $f_n=\frac{1}{\sigma_n\sqrt{n\log n}}
\sum_{k=0}^{n-1} {\bf 1}_{[k,k+1]}^{\otimes q}
\in\HH^{\odot q}$. According to Lemma \ref{Stein-lemma},  we have,
for all $k\geq 1$ and $r=1,\ldots,q-1$,
$\|f_k\otimes_r f_k\|_{\HH^{\otimes (2q-2r)}}\leq C/\sqrt{\log k}$.
Hence
\begin{eqnarray*}
\sum_{n\geq 2}\frac{1}{n\log^2 n}\sum_{k=1}^n \frac{1}{k}\|f_k\otimes_r f_k\|_{\HH^{\otimes (2q-2r)}}
&\leq&
C\sum_{n\geq 2}\frac{1}{n\log^2 n}\sum_{k=1}^n \frac{1}{k\sqrt{\log k}}\\
&\leq&
C\sum_{n\geq 2}\frac{1}{n\log^{3/2} n}
<\infty,
\end{eqnarray*}
that is assumption $(A_1)$ is verified.
Concerning $(A_2)$, note that
\[
\langle f_k,f_l\rangle_{\HH^{\otimes q}} =\frac{1}{\sigma_k\sigma_l\,\sqrt{k\log k}\sqrt{ l\log l}}
\sum_{i=0}^{k-1}\sum_{j=0}^{l-1}\rho(j-i)^q.
\]
Since $\sigma_k\rightarrow\sigma_\infty>0$ as $k\to\infty$ (see Lemma \ref{sigma}), we have, 
for $l\geq k\geq 1$,
\begin{eqnarray*}
\big|\langle f_k,f_l\rangle_{\HH^{\otimes q}}\big|&\leq&\frac{C}{\sqrt{k\log k}\sqrt{l\log l}}
\sum_{i=0}^{k-1}\sum_{j=0}^{l-1}\big|\rho(j-i)\big|^q\\
&=&
\frac{C}{\sqrt{k\log k}\sqrt{l\log l}}
\sum_{i=0}^{k-1}\sum_{r=-i}^{l-1-i}\big|\rho(r)\big|^q\\
&\leq&C\frac{\sqrt{k}}{\sqrt{\log k}\sqrt{l\log l}}\,\sum_{r=-l}^{l}\big|\rho(r)\big|^q 
\leq C\sqrt{\frac{k\log l}{l\log k}}.
\end{eqnarray*}
The last inequality follows from the fact that $\sum_{r=-l}^l\big|\rho(r)\big|^q\leq C\log l$ 
since, as $|r|\to\infty$, 
\[
\rho(r)\sim(1-\frac{1}{q})(1-\frac1{2q})|r|^{-1/q}.
\]
Consequently, assumption $(A_2)$ is verified because
\begin{eqnarray*}
\sum_{n\geq 2}\frac{1}{n\log^3 n}\sum_{k,l=2}^{n} \frac{\big|\langle f_k,f_l\rangle_{\HH^{\otimes q}}\big|}{kl}
&\leq&
2\sum_{n\geq 2}\frac{1}{n\log^3 n}\sum_{l=2}^{n}\sum_{k=2}^l \frac{\big|\langle f_k,f_l\rangle_{\HH^{\otimes q}}\big|}{kl}\\
&\leq&C\sum_{n\geq 2}\frac{1}{n\log^3 n}\sum_{l=2}^n \frac{\sqrt{\log l}}{l^{3/2}}\sum_{k=2}^l \frac1{\sqrt{k\log k}}\\
&\leq&C\sum_{n\geq 2}\frac{1}{n\log^3 n}\sum_{l=2}^n \frac1l   \leq 
C\sum_{n\geq 2}\frac{1}{n\log^2 n}<\infty.
\end{eqnarray*}
\fin

\section{Partial sums of Hermite polynomials of increments of fBm: the non-Gaussian case}\label{s:hermite}
Fix $q\geq 2$. In the previous section, we saw that the limit distribution
of $V_n=\sum_{k=0}^{n-1} H_q(B^H_{k+1}-B^H_{k})$, adequably normalized, is Gaussian
when $H\leq 1-1/(2q)$. We consider here
the case $H>1-1/(2q)$.
In contrast to (\ref{eq:Breuer_Major1})-(\ref{eq:Breuer_Major2}), we have
\begin{equation}
\label{eq:Breuer_Major3}
G_n:=n^{q(1-H)-1}V_n 
\,\overset{{\rm Law}}{\longrightarrow}\,
G_\infty.
\end{equation}
The law of $G_\infty$ is called the ``Hermite distribution''. 
A short proof of (\ref{eq:Breuer_Major3}) is given in Proposition \ref{lm-hermite-yahoo} below. 
It is based on the fact that, for {\it fixed} $n$, 
$Z_n$ defined in (\ref{sn}) below and $G_n$
share the same law, 
because of the self-similarity property of fractional Brownian motion. 
\begin{prop}\label{lm-hermite-yahoo}
Fix $q\geq 2$ and $H>1-1/(2q)$, and define $Z_n$
by
\begin{equation}\label{sn}
Z_n=n^{q(1-H)-1}\sum_{k=0}^{n-1} H_q\big(n^H(B^H_{(k+1)/n}-B^H_{k/n})\big),\quad n\geq 1.
\end{equation}
Then, as $n\to\infty$,
$\{Z_n\}$ converges almost surely and in $L^2(\Omega)$ to a limit denoted by $Z_\infty$.
\end{prop}
{\bf Proof}. 
Let us first prove the convergence in $L^2(\Omega)$.
For $n,m\geq 1$, we have
\[
E[Z_nZ_m]=q!(nm)^{q-1}\sum_{k=0}^{n-1}\sum_{l=0}^{m-1}\left(E\big[\big(B^H_{(k+1)/n}-B^H_{k/n}\big)
\big(B^H_{(l+1)/m}-B^H_{l/m}\big)\big]\right)^q.
\]
On the other hand, since $H>1/2$, we have, for all $s,t\geq 0$,
\[
E[B^H_sB^H_t]=H(2H-1)\int_0^t du\int_0^s dv |u-v|^{2H-2}.
\]
Hence
\[
E[Z_nZ_m]=q!H^q(2H-1)^q\times\frac1{nm}
\sum_{k=0}^{n-1}\sum_{l=0}^{m-1}\left(nm\int_{k/n}^{(k+1)/n} du \int_{l/m}^{(l+1)/m} dv 
|v-u|^{2H-2}\right)^q.
\]
Therefore, as $n,m\to\infty$, we have, 
\[
E[Z_nZ_m]\to
q!H^q(2H-1)^q
\int_{[0,1]^2}|u-v|^{(2H-2)q}dudv,
\]
and the limit is finite since $H>1-1/(2q)$.
In other words, the sequence $\{Z_n\}$ is Cauchy in $L^2(\Omega)$, 
and hence converges in $L^2(\Omega)$ to some $Z_\infty$. 

Let us now prove that $\{Z_n\}$ converges also almost surely. 
Observe first that, since $Z_n$ belongs to the $q$th chaos of $B^H$ for all $n$ and because $\{Z_n\}$
converges in $L^2(\Omega)$ to $Z_\infty$, we have that $Z_\infty$ also belongs to the $q$th chaos of $B^H$.
In \cite[Proposition 3.1]{BN}, it is shown that $E[|Z_n-Z_\infty|^2]\leq Cn^{2q-1-2qH}$, for some 
positive constant $C$ not depending on $n$. 
Inside a fixed chaos, all the $L^p(\Omega)$-norms are equivalent
(see e.g. \cite[Theorem 5.10]{Janson}). Hence, for any $p>2$,
we have $E[|Z_n-Z_\infty|^p]\leq Cn^{p(q-1/2-qH)}$.
Since $H>1-1/(2q)$, there exists $p>2$ large enough such that
$(q-1/2-qH)p<-1$. Consequently 
\[
\sum_{n\geq 1}E[|Z_n-Z_\infty|^p]<\infty
\]
leading, for all $\e>0$, to 
\[
\sum_{n\geq 1}P[|Z_n-Z_\infty|>\e]<\infty.
\]
Therefore, we deduce from the Borel-Cantelli lemma
that $\{Z_n\}$ converges 
almost surely to $Z_\infty$. 
\fin

We now face some difficulties. First,
since the limit of $G_n$ in (\ref{eq:Breuer_Major3}) is not Gaussian, we cannot apply our general
criterion Theorem \ref{main} to obtain an ASCLT.
To modify adequably the criterion, we would need a version of Proposition \ref{noupec}
for random variables with a Hermite distribution, a result which is not presently available.
Thus an ASCLT associated to the convergence in law (\ref{eq:Breuer_Major3}) falls outside the
scope of this paper. We can nevertheless make a number of observations. 
First, changing the nature of the random variables without changing their law has no impact
on CLTs as in (\ref{eq:Breuer_Major3}), but may have a great impact on an ASCLT. To see this,
observe that for each fixed $n$, the ASCLT involves not only the distribution of $G_n$, but also
that of $(G_1,\ldots,G_n)$. 
Consider, moreover, the following example. Let $\{G_n\}$ be a sequence of random variables converging in law
to a limit $G_\infty$. According to a theorem of Skorohod, there is a sequence $\{G_n^*\}$ such that
for any fixed $n$, $G_n^*\overset{\rm Law}{=}G_n$ and such that $\{G_n^*\}$ converges almost surely
as $n\to\infty$ to a random variable $G^*_\infty$ with $G_\infty^*\overset{\rm Law}{=}G_\infty$.
Then, for any bounded continuous function $\varphi:\R\to\R$, we have
$\varphi(G_n^*)\overset{\rm a.s.}{\longrightarrow}\varphi(G_\infty^*)$ and by 
a classical theorem of Hardy (see \cite{chmi} p.35), as $n\to\infty$,
\[
\frac1{\log n}\sum_{k=1}^n \frac1k \varphi(G_k^*)\overset{\rm a.s.}{\longrightarrow}\varphi(G_\infty^*).
\] 
This limit is different from $E[\varphi(G_\infty^*)]$ (or equivalently 
$E[\varphi(G_\infty)]$),
that is, different from the limit if one had an ASCLT. Thus, knowing the law of $G_n$, for a fixed $n$,
does not allow to determine whether an ASCLT holds or not.

\begin{rem}
{\rm
In view of Proposition \ref{lm-hermite-yahoo}, the Skorohod
version of
\[
G_n=n^{q(1-H)-1}\sum_{k=0}^{n-1}H_q(B^H_{k+1}-B^H_k)
\]
is
\[
G_n^*=Z_n=n^{q(1-H)-1}\sum_{k=0}^{n-1}H_q\big(
n^H(B^H_{(k+1)/n}-B^H_{k/n})\big),
\]
since $G_n^*\overset{\rm law}{=}G_n$ and 
$G_n^*$ converges almost surely.
}
\end{rem}

Hence, in the case of Hermite distributions, by suitably modifying the argument of the Hermite polynomial
$H_q$ in a way which does not change the limit in law, namely by considering $Z_n$ in (\ref{sn})
instead of $G_n$ in (\ref{eq:Breuer_Major3}), we obtain the almost sure convergence
\[
\frac1{\log n}\sum_{k=1}^n\frac1k \varphi(Z_k)\to 
\varphi(Z_\infty).
\]
Note that the limit is different from the limit
expected under an ASCLT, namely $E[\varphi(Z_\infty)]$. Recall indeed that $Z_\infty$ is a non-constant
random variable with a Hermite distribution
(Dobrushin and Major \cite{DoMa}, Taqqu \cite{T}) and therefore
one has $E[\varphi(Z_\infty)]\neq \varphi(Z_\infty)$ in general.\\

\bigskip

\noindent
{\bf Acknowledgments}. This paper originates from the conference ``Limit theorems and applications'',
University Paris I Panth\'eon-Sorbonne,
January 14-16, 2008, that the three authors were attending. We warmly thank J.-M. Bardet and C. A. Tudor
for their invitation and generous support. Also, 
I. Nourdin would like to thank M. S. Taqqu
for his hospitality during his stay at Boston University in March 2009, where part of this research was carried out.


\begin{thebibliography}{99}

\bibitem{BerkesCsaki}
I. Berkes and E. Cs\'aki (2001).
A universal result in almost sure central limit theory.
{\it Stoch. Proc. Appl.} {\bf 94}, no. 1, 105-134.


\bibitem{BerkesHorvath}
I. Berkes and L. Horv{\'a}th (1999).
Limit theorems for logarithmic averages of fractional Brownian motions.
{\it J. Theoret. Probab.} {\rm 12}, no. 4, 985--1009.


\bibitem{BN}
J.-C. Breton and I. Nourdin (2008).
Error bounds on the non-normal approximation of Hermite power variations of fractional Brownian motion.
{\it Electron. Comm. Probab.} {\bf 13}, 482-493.


\bibitem{BM}
P. Breuer and P. Major (1983). 
Central limit theorems for nonlinear functionals of Gaussian fields. 
{\it J. Multivariate Anal.} {\bf 13}, no. 3, 425-441.

\bibitem{chmi}
K. Chandrasekharan and S. Minakshisundaram (1952).
{\it Typical means}.
Oxford University Press, Oxford.

\bibitem{DoMa} 
R. L. Dobrushin and P. Major (1979).
Non-central limit theorems for nonlinear functionals of Gaussian fields.
{\it Z. Wahrsch. verw. Gebiete}, no. 50, 27-52. 

\bibitem{GS} 
L. Giraitis and D. Surgailis (1985). 
CLT and other limit theorems for functionals of Gaussian processes. 
{\it Z. Wahrsch. verw. Gebiete}, no. 70, 191-212. 

\bibitem{Brosamler}
G. A. Brosamler (1988).
An almost everywhere central limit theorem.
{\it Math. Proc. Cambridge Philos. Soc.} {\bf 104}, no. 3, 561-574.

\bibitem{these}
K. Gonchigdanzan (2001).
Almost Sure Central Limit Theorems.
PhD thesis, University of Cincinnati.
Available online at {\tt www.ohiolink.edu/etd/view.cgi?acc\_\,num=ucin990028192}.

\bibitem{IL}
I. A. Ibragimov and M. A. Lifshits (2000).
On limit theorems of ``almost sure'' type.
{\it Theory Probab. Appl.} {\bf 44}, no. 2, 254-272. 

\bibitem{IbragimovLifshits2}
I. A. Ibragimov and M. A. Lifshits (1998).
On the convergence of generalized moments in almost sure central limit theorem.
{\it Statist. Probab. Lett.} {\bf 40}, no. 4, 343-351.

\bibitem{Janson}
S. Janson (1997).
{\it Gaussian Hilbert Spaces}.
Cambridge University Press.


\bibitem{LaceyPhillip}
M.T. Lacey and W. Philipp (1990).
A note on the almost sure central limit theorem.
{\it Statist. Probab. Letters} {\bf 9}, 201-205.

\bibitem{Levy}
P. L\'evy (1937).
{\it Th\'eorie de l'addition des variables al\'eatoires}.
Gauthiers-Villars. 

\bibitem{NP07}
I. Nourdin and G. Peccati (2007). 
Stein's method on Wiener chaos.
{\it Probab. Theory Related Fields}, to appear.

\bibitem{NP} 
D. Nualart and G. Peccati (2005). 
Central limit theorems for sequences of multiple stochastic integrals. 
{\it Ann. Probab.} {\bf 33}, no. 1, 177-193.


\bibitem{Nbook}
D. Nualart (2006). 
{\it The Malliavin calculus and related topics.}
Springer-Verlag, Berlin, 2nd edition.

\bibitem{SamorTaqqu}
G. Samorodnitsky and M. S. Taqqu (1994).
{\it Stable non-Gaussian random processes}.
Chapman and Hall, New York. 

\bibitem{Schatte}
P. Schatte (1988).
On strong versions of the central limit theorem.
{\it Math. Nachr.} {\bf 137}, 249-256.

\bibitem{T} 
M. S. Taqqu (1979).
Convergence of integrated processes of arbitrary Hermite rank.
{\it Z. Wahrsch. verw. Gebiete} {\bf 50}, 53-83. 


\end{thebibliography}
\end{document}